\documentclass[12pt]{article}
\usepackage[final]{epsfig}
\usepackage{graphics}
\usepackage{amsmath}
\usepackage{amsfonts}
\usepackage{latexsym}
\usepackage{amssymb}
\usepackage{graphicx}
\usepackage{url}
\usepackage{epstopdf}

\newtheorem{lemma}{Lemma}[section]

\newtheorem{theorem}{Theorem}

\newtheorem{corollary}[theorem]{Corollary}

\begin{document}
\newcommand{\eps}{{\varepsilon}}
\newcommand{\proofend}{$\Box$\bigskip}
\newcommand{\C}{{\mathbb C}}
\newcommand{\Q}{{\mathbb Q}}
\newcommand{\R}{{\mathbb R}}
\newcommand{\Z}{{\mathbb Z}}
\newcommand{\RP}{{\mathbb {RP}}}
\newcommand{\CP}{{\mathbb {CP}}}
\newcommand{\Tr}{\rm Tr}
\def\proof{\paragraph{Proof.}}
\newcommand{\SL}{\operatorname{SL}}
\newcommand{\PGL}{\operatorname{PSL}}

\title{Descartes Circle Theorem, Steiner Porism, and Spherical Designs}

\author{Richard Evan Schwartz\footnote{
Department of Mathematics, Brown University, Providence, RI 02912;
res@math.brown.edu} \and 
Serge Tabachnikov\footnote{
Department of Mathematics,
Penn State University,
University Park, PA 16802;
tabachni@math.psu.edu}
}

\date{}
\maketitle

\section{Introduction} \label{sect:intro}

The Descartes Circle Theorem has been popular lately because it underpins the geometry and arithmetic of Apollonian packings, a subject of great current interest; see, e.g., surveys \cite{Fu,Kon,LMW,Sar}. In this article we revisit this classic result, along with another old theorem on circle packing, the Steiner porism, and relate these topics to spherical designs. 

The Descartes Circle Theorem concerns cooriented circles. A coorintation of a circle is a choice of one of the two components of its complement. One can think of a normal vector field along the circle pointing toward this component. If a circle has radius $r$, its curvature (bend) is given by the formula $b=\pm 1/r$, where the sign is positive if the normal vector is inward and negative otherwise. By tangent cooriented circles we mean tangent circles whose coorientations are opposite; that is, the normal vectors have opposite directions at the point of tangency. 

As usual in inversive geometry, we think of straight lines as circles of infinite radius and zero curvature. Our definitions also extend to spheres in higher dimensional Euclidean spaces; they include spheres of infinite radius, that is, hyperplanes. 

The Descartes Circle Theorem relates the signed curvatures $a,b_1,b_2,b_3$ of four mutually tangent cooriented circles:
\begin{equation} \label{Desrule}
(a+b_1+b_2+b_3)^2-2(a^2+b_1^2+b_2^2+b_3^2)=0.
\end{equation} 
Considered as a quadratic equation in $a$, equation (\ref{Desrule}) has two roots
\begin{equation*}
a_{1,2}=(b_1+b_2+b_3) \pm 2\sqrt D, \hskip 30 pt
D=b_1b_2+b_2b_3+b_3b_1,
\end{equation*}
and hence (Vieta formulas)
\begin{equation} \label{Soddy}
a_1+a_2=2(b_1+b_2+b_3),\ \ a_1 a_2 = b_1^2+b_2^2+b_3^2-2(b_1b_2+b_2b_3+b_3b_1).
\end{equation}

\begin{figure}[hbtp]
\centering
\includegraphics[height=2.2in]{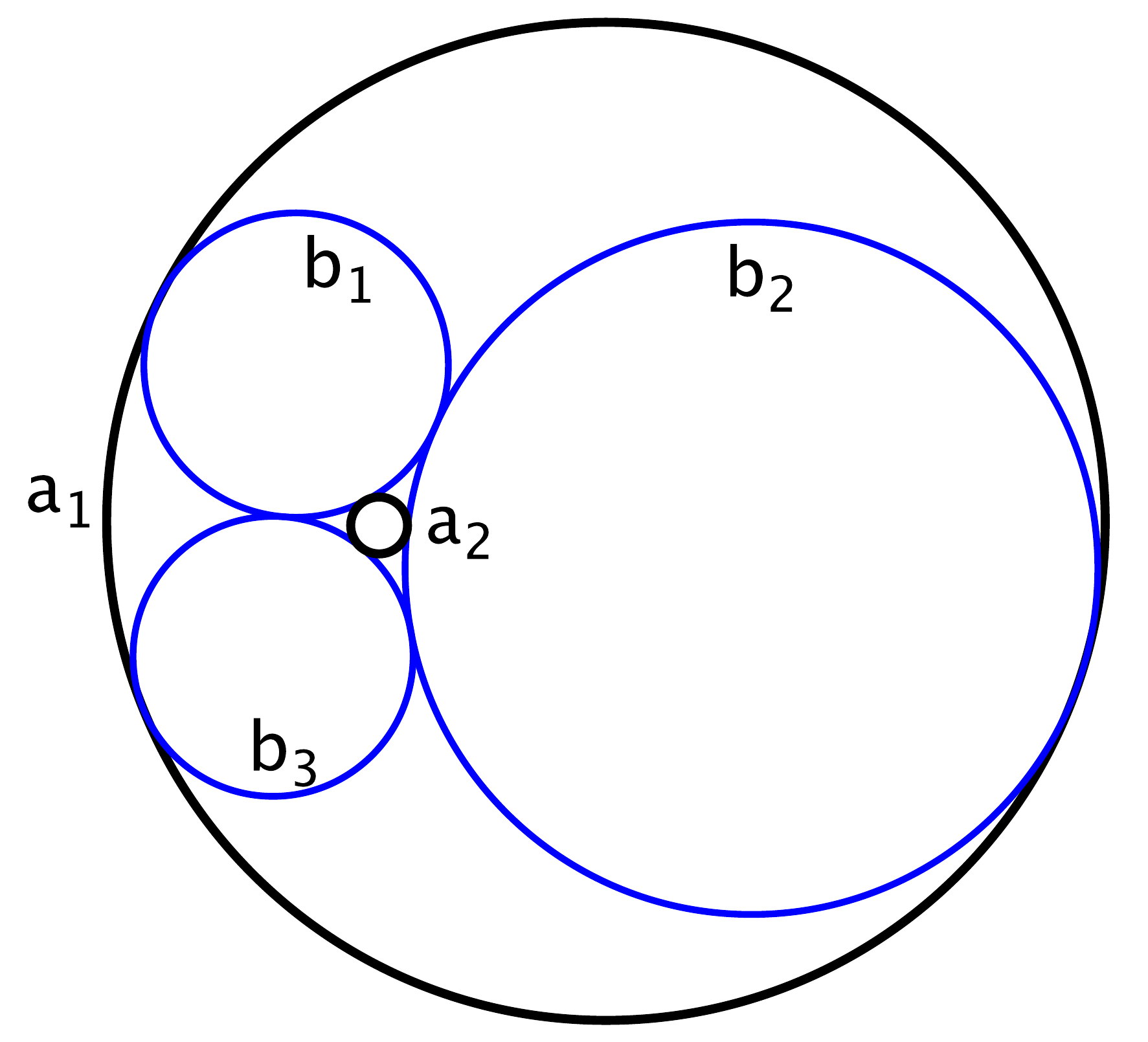}\qquad 
\includegraphics[height=2.2in]{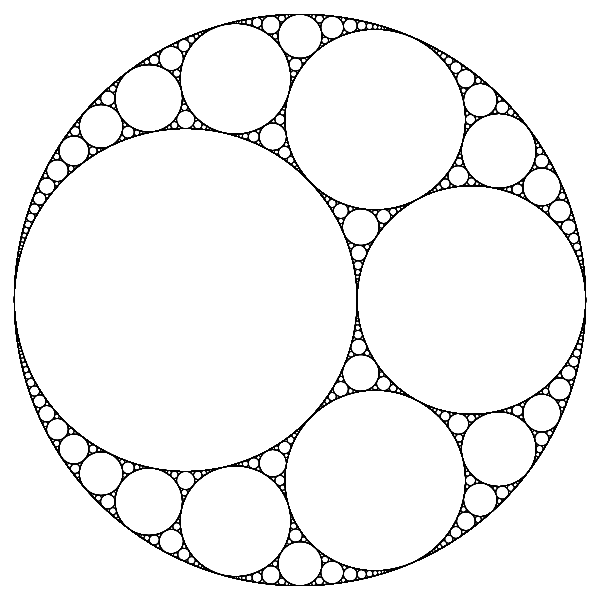}
\caption{Three mutually tangent circles inscribed into an annuls and an Apollonian gasket.}
\label{five}
\end{figure}

Figure \ref{five} illustrates the situation: the curvatures of the inner-most and outer-most circles are $a_{1,2}$, and the curvatures of the three circles in-between are $b_{1,2,3}$. 
Thus, starting with a Descartes configuration of four mutually tangent circles, one can select one of the circles and replace it by another circle, tangent to the same three circles (replace $a_1$ by $a_2$). Continuing this process results in an Apollonian packing of circles (Apollonian gasket).

The  connection to arithmetic is that if $a_1,b_1,b_2,b_3$ are all integers then, by the first equation in (\ref{Soddy}),
so is $a_2$.  Thus, if four mutually tangent circles in an Apollonian packing have integer curvatures,
then all the circle do.  This fact opens the door to many deep and subtle works on the arithmetic of the Apollonian packing.

One can rewrite formulas (\ref{Soddy}) as 
\begin{equation} \label{moments1}
b_1+b_2+b_3=\frac{a_1+a_2}{2},\ \ b_1^2+b_2^2+b_3^2=\frac{6a_1a_2-a_1^2-a_2^2}{4}.
\end{equation}
We note that a result like this could not work for the third moments of curvature
because it would force $b_1,b_2,b_3$ to be independent of the choices of the respective circles, clearly, not the case.

In fact, there is one degree of freedom in choosing the triple of mutually tangent circles contained in the annulus between the inner- and outer-most circles and tangent to them.  This is the a particular case of the Steiner porism.

\begin{figure}[hbtp]
\centering
\includegraphics[height=2.5in]{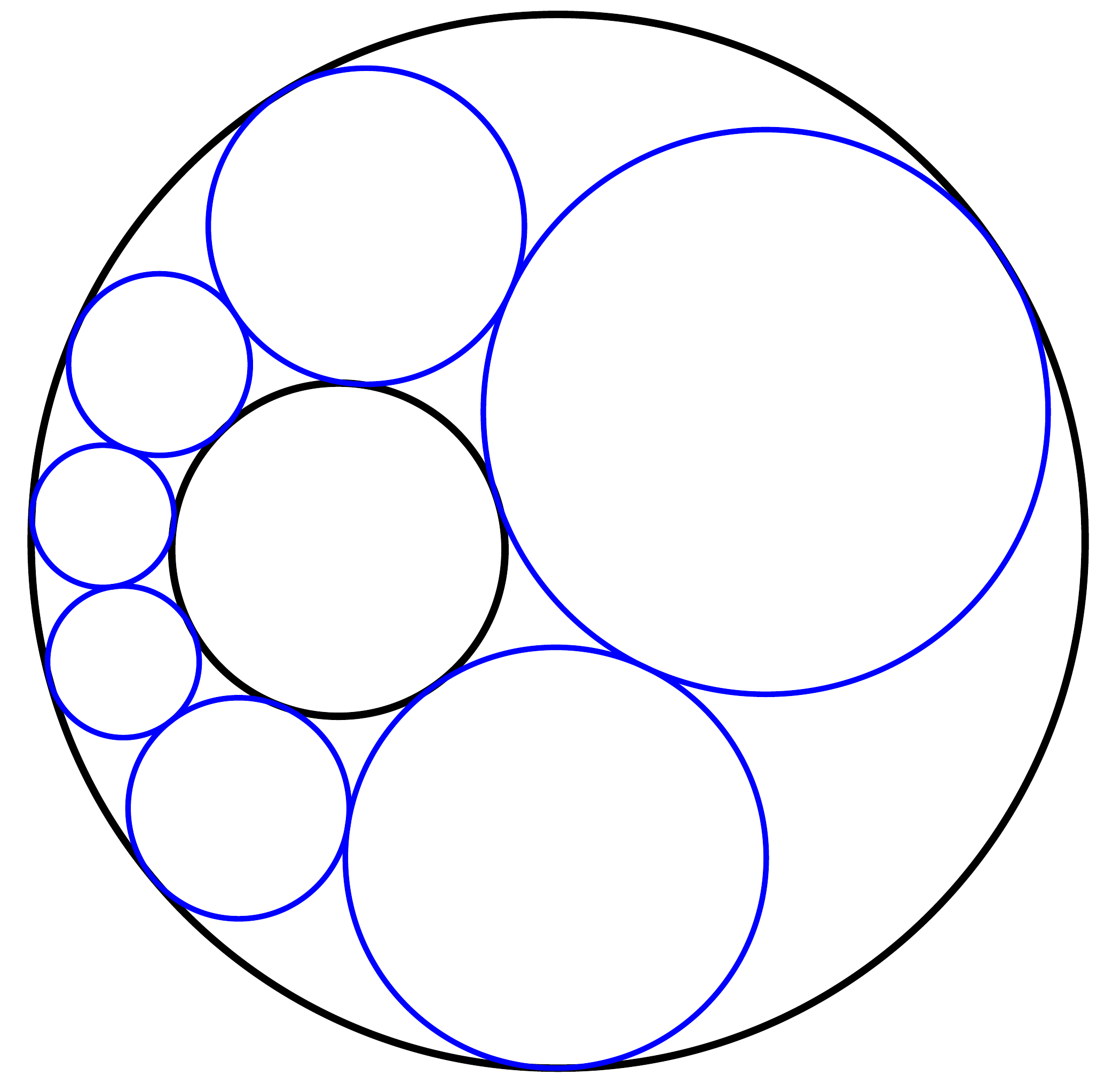}
\caption{A Steiner chain of length 7.}
\label{seven}
\end{figure}

Given two nested circles, inscribe a circle in the annulus bounded by them, and construct a {\it Steiner chain of circles},  tangent to each other in a cyclic pattern, see Figure \ref{seven}. The claim is that if this chain closes up after $k$ steps (perhaps making several turns around the annulus), then the same will happen for any choice of the initial circle. In other words, one can rotate a Steiner chain around the annulus, like a ball bearing. 
Let us call the inner-most and outer-most circles the {\it parent circles}.

The Steiner porism is easy to prove: there exists an inversion that takes the parent circles to concentric circles (see, e.g., \cite{Cox}, section 6.5). Since an inversion takes circles to circles and preserves tangency, the theorem reduces to the case of concentric circles, where it follows from  rotational symmetry. In other words, there is a group of conformal transformations, isomorphic to $SO(2)$, preserving the parent circles, and sending  Steiner chains to Steiner chains of the same length (and the same  number of turns about the annulus). 

The Steiner porism is the conformal geometry analog of the famous Poncelet porism from projective geometry (see, e.g., \cite{Fla}).  The Poncelet porism is a deeper result however, because the proof does not just boil down to an application of projective symmetry.\footnote {See \cite{Bar} for an analog of the Full Poncelet Theorem for Steiner chains.} 

Suppose that a pair of parent circles admits a Steiner chain of length $k \ge 3$. Then there is a 1-parameter family of Steiner chains of length $k$ with the same parent circles. What do the $k$ circles in these chains have in common? The next result is a generalization of formulas (\ref{moments1}).

\begin{theorem} \label{Chain1}
For every $m=1,2,\ldots,k-1$, the sum
$$
I_m=\sum_{j=1}^k  b_j^m 
$$
remains constant in the 1-parameter family of Steiner chains of length $k$. 
\end{theorem}

One can prove that $I_m$ is a symmetric homogeneous polynomial of degree $m$ in the curvatures of the parent circles (symmetry follows from the fact that there exists a conformal involution that interchanges the parent circles and sends Steiner chains to Steiner chains). We do not dwell on this addition to Theorem \ref{Chain1} here.

One may reformulate Theorem \ref{Chain1} as follows. Associate with a Steiner chain of length $k$ with fixed parent circles  the polynomial $P(x)=(x-b_1)(x-b_2)\cdots (x-b_k)$. One obtains a 1-parameter family of polynomials of degree $k$. 

\begin{corollary} \label{poly}
The polynomials $P(x)$ in this 1-parameter family differ only by constants. That is, the first $k-1$ symmetric polynomials of the curvatures $b_1,\ldots,b_k$ remain the same, and only the product $b_1\cdots b_k$ varies in the 1-parameter family.
\end{corollary}

\proof
The coefficients of $P(x)$ are the elementary symmetric functions of $b_j$. One has Newton identities between elementary symmetric functions $\sigma_m$ and powers of sums $I_m$ (as in Theorem \ref{Chain1}):
$$
\sigma_1=I_1,\ 2\sigma_2=\sigma_1I_1-I_2,\ 3\sigma_3=\sigma_2I_1-\sigma_1I_2+I_3,\ 4\sigma_4=\sigma_3I_1-\sigma_2 I_2+\sigma_1I_3-I_4,
$$
and so on (see, e.g., \cite{Me}). By Theorem \ref{Chain1}, $I_1,\ldots, I_{k-1}$ remain constant in the 1-parameter family of Steiner chains, and hence so do $\sigma_1,\ldots, \sigma_{k-1}$.  
\proofend

The Descartes Circle Theorem has a recent generalization, discovered in \cite{LMW}. 
Let $w,z_1,z_2,z_3 \in \C$ be the centers of the four mutually tangent circles,  thought of as complex numbers, and let 
their respective curvatures be $a,b_1,b_2.b_3$. Then one has the Complex Descartes Theorem
\begin{equation*} \label{comDesrule}
(aw+b_1z_1+b_2z_2+b_3z_3)^2-2(a^2w^2+b_1^2z_1^2+b_2^2z_2^2+b_3^2z_3^2)=0,
\end{equation*}
an analog of (\ref{Desrule}). 
As before, this implies that the two complex numbers, 
\begin{equation} \label{moments2}
b_1z_1+b_2z_2+b_3z_3 \ \ {\rm and}\ \ b_1^2z_1^2+b_2^2z_2^2+b_3^2z_3^2
\end{equation}
remain constant in the 1-parameter family of Steiner chains of length 3 with fixed parent circles. 

Relations (\ref{moments2}) remain valid if one parallel translates the configuration of circles, that is, adds an arbitrary complex number $u$ to  each $z_i$. Taking the coefficients of the respective powers of $u$ implies that the first and the second moments of the curvatures $b_1,b_2,b_3$ remain constant in the 1-parameter family of Steiner chains -- as we already know, see (\ref{moments1}), and  yields another conserved quantity
$$
b_1^2z_1+b_2^2z_2+b_3^2z_3.
$$

We generalize to Steiner chains of length $k$.
Let $z_1,\ldots,z_k \in \C$ be the centers of the circles in such a chain (these centers lie on a conic whose foci are the centers of the parent circles).

\begin{theorem} \label{Chain2}
For all $0\le n \le m \le k-1$, the sum
$$
J_{m,n}=\sum_{j=1}^k  b_j^m z_j^n 
$$
remains constant in the 1-parameter family of Steiner chains of length $k$. 
\end{theorem}

Since $I_m=J_{m,0}$, Theorem \ref{Chain2} implies Theorem \ref{Chain1}. Since the circles that form a Steiner chain with fixed parent circles have only one degree of freedom, the quantities  $J_{m,n}$ are not independent and satisfy  algebraic relations. As before, it suffices to establish Theorem \ref{Chain2} for $n=m$ and apply parallel translation to obtain its general statement.  

\section{Spherical designs} \label{design}

Let $S^d \subset \R^{d+1}$ be the unit  origin centered sphere. A  collection of points $X=\{\bar x_1,\ldots,\bar x_n\} \subset S^d$ is called a {\it spherical $M$-design} if, for every polynomial function on $\R^{d+1}$ of degree not greater than $M$, its average value over $S^d$ equals its average over $X$ (the average is taken with respect to the standard, rotation-invariant measure on the sphere). 


Here are some examples:
\begin{enumerate}
\item Take $d=1$. If $X$ consists of $M+1$ evenly
spaced points on $S^1$, then $X$ is a spherical
$M$-design.  This derives from the fact that the
sums of the $m$th powers of the $(M+1)$th roots of unity
are $0$ for $m=1,...,M$.
\item Take $d=2$. The vertices of the
tetrahedron, cube, octahedron, dodecahedron,
icosahedron respectively are spherical $M$-designs
for  $M=2,3,3,5,5$.
\item Take $d=3$. The set of vertices of the $24$-cell,
the $120$-cell, and the $600$-cell respectively
are $M$-designs for $M=5,11,11$.
\item Take $d=7$.  The vertices of the $E_8$ cell form a $7$-design.
\item Take $d=23$. The vertices of the Leech cell form an $11$-design.
\end{enumerate}
For more information, see the survey \cite{Ba}.

\begin{figure}[hbtp]
\centering
\includegraphics[height=2in]{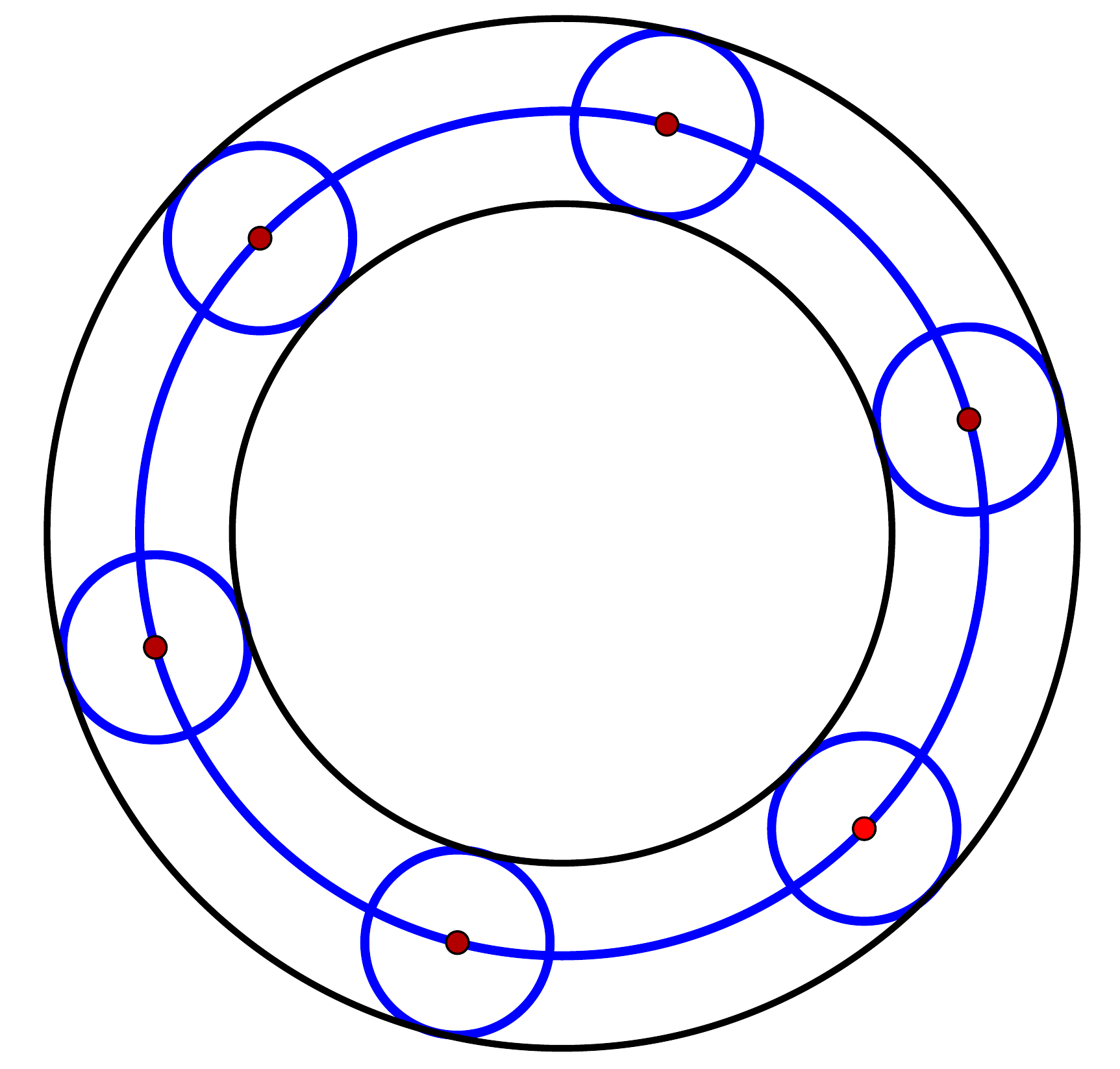}\qquad 
\includegraphics[height=2in]{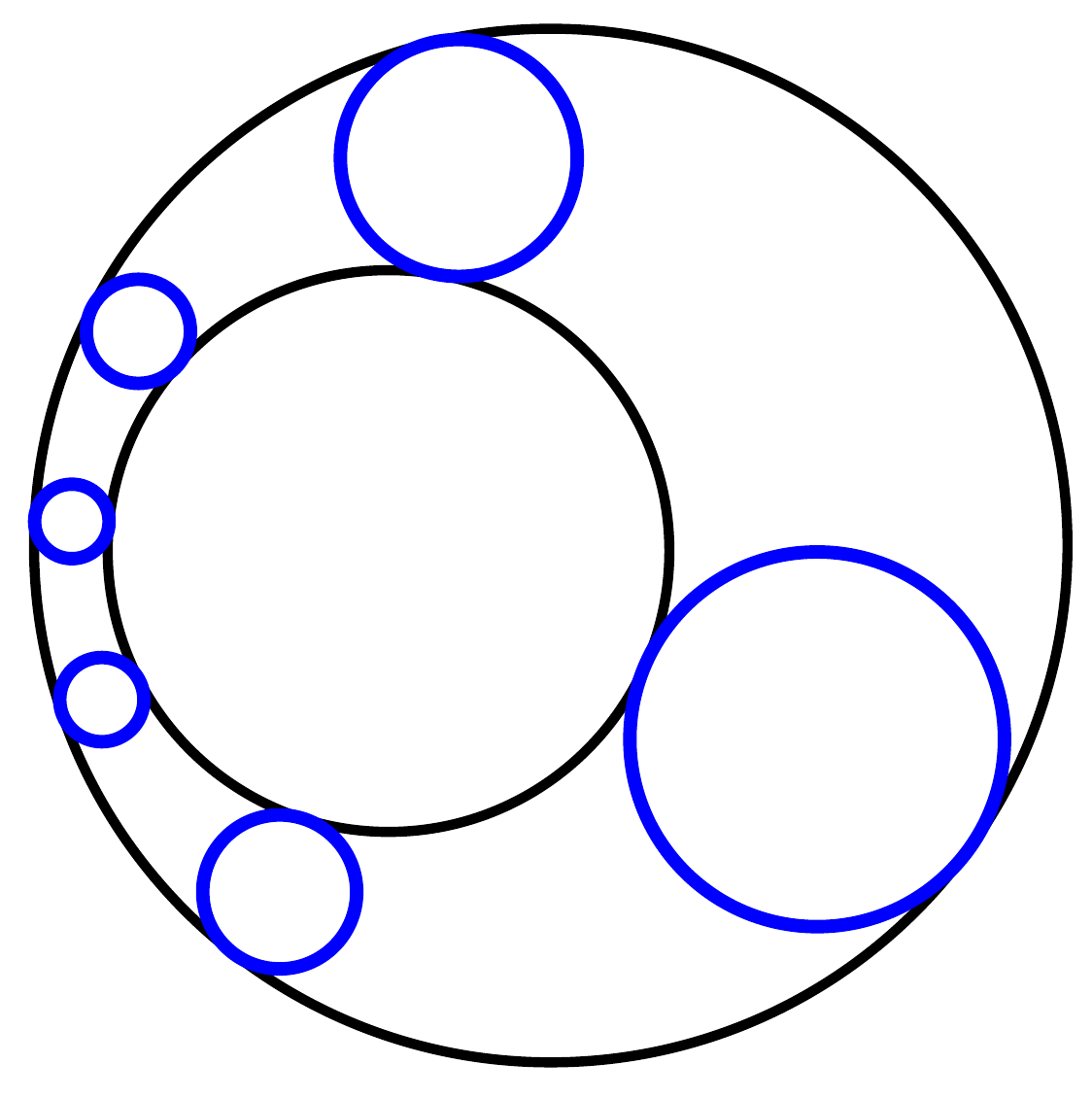}
\caption{A ball bearing with gaps and its conformal image.}
\label{gaps}
\end{figure}

Let $S^d \subset \R^{d+1}$ be the unit sphere, and $X=\{\bar x_1,\ldots,\bar x_n\}$ be a {spherical $M$-design}. Fix a positive number $\bar b$ and consider the spheres of curvature $\bar b$, centered at points $\bar x_i,\ i=1,\ldots,n$. For example, if $d=1$, one has  $M=n-1$, and we have $n$ equal circles whose centers are evenly spaced on the circle (a ball bearing with gaps, see Figure \ref{gaps}).

Let $Y=\{\bar y_1,\ldots, \bar y_q\} \subset S^d$ be another spherical $M$-design. Again consider the spheres of the same curvature $\bar b$, centered at points $\bar y_j,\ j=1,\ldots,q$. For example, if $q=n$, then $Y$ could be obtained from $X$ by an isometry of the sphere (rotation of the circle, for $d=1$). 

Let $\Psi$ be conformal transformation. Apply $\Psi$ to the two sets of spheres. Let $(x_i, b_i)$ be the centers and the curvatures for the first set of spheres, and $(y_j, c_j)$ be the centers and the curvatures for the second set. 

Let $F$ be a polynomial on $\R^{d+1}$ of degree not greater than $M$. The next result generalizes Theorem \ref{Chain2}. 

\begin{theorem} \label{moments}
One has 
$$
\frac{1}{n} \sum_{i=1}^n F(b_i  x_i) = \frac{1}{q} \sum_{j=1}^q F( c_j  y_j).
$$
\end{theorem}

For the proof, we shall need two lemmas, both can be found in the literature on inversive geometry (see, e.g., \cite{Ped,Wil}).

\begin{lemma} \label{struct}
Every conformal transformation of $\R^{d+1}$ is either a similarity or a composition of similarity with the inversion in a unit sphere.
\end{lemma}

\proof
A conformal transformation that sends infinity to infinity is a similarity. If $\Psi$ sends $\infty$ to point $C$, let $I$ be the inversion in the unit sphere centered at $C$. Then $I(C)=\infty$, and $S:= I\circ \Psi$ sends infinity to infinity.  Hence $\Psi=I\circ S$, where $S$ is a similarity.
\proofend

The case of Theorem \ref{moments} when $\Psi$ is a similarity is trivial since similarities preserve the space of polynomials of degree $\le M$. To tackle the general case, we need to know how inversion affects the curvatures and centers of spheres.

Consider a sphere centered at point $\bar x \in \R^{d+1}$ with curvature $\bar b$. Invert this sphere in the unit sphere centered at point $C$ and let $x,b$ be the center and the curvature of the image sphere. The next lemma is a result of a computation that we omit.

\begin{lemma} \label{invform}
One has:
$$
b=\bar b |\bar x-C|^2 - \frac{1}{\bar b},\ x = C + \frac{\bar x - C}{|\bar x-C|^2 - \frac{1}{\bar b^2}},\ 
bx = \bar b \bar x + \bar bC \left( |\bar x-C|^2 - \frac{1}{\bar b^2} -1 \right).
$$
\end{lemma}

Now we can prove Theorem \ref{moments}.

\proof
We have $\Psi=I\circ S$, where $I$ is the inversion in the unit sphere centered at point $C$. 

The original configuration of $n$ congruent spheres centered at the points of a spherical $M$-design is taken by 
the similarity $S$ to a configuration of $n$ congruent spheres of the same curvature $\bar b$ whose centers $\bar x_i,\ i=1,\ldots,n$, lie on a sphere of some radius $r$, centered at some point $A$. Thus $\bar x_i = A + r v_i$ where $v_1,\ldots,v_n$ is a spherical $M$-design. 

If $\bar x = A + r v$ with $v \in S^d$, then
$$
|\bar x - C|^2 - \frac{1}{\bar b^2} -1= |A-C|^2 + r^2 - \frac{1}{\bar b^2} + 2r (A-C) \cdot v,
$$
the sum of a constant and a linear function, that is, a polynomial of degree one on $S^d$ (depending on $A,C,r,\bar b$). 

Now the last equation of Lemma \ref{invform} implies that, for every $k=1,\ldots,d+1$, the $k$th component of the vector $b_i x_i$ is given by a polynomial of degree one $\alpha_k + \beta_k(v_i)$, where $\alpha_k$ is a constant and $\beta_k$ is a linear function on $S^d$ (depending on $A,C,r,\bar b$, and $k$). 

Therefore the function $F$ is a polynomial of degree at most $M$ on $S^d$, and its average over the set $v_1,\ldots,v_n$ equals its average over the sphere. The same applies to the second spherical $M$-design consisting of $q$ points, and the result follows.  
\proofend



One also has an analog of Theorem \ref{Chain1}. 

\begin{corollary} \label{desmom}
Under the assumption of Theorem \ref{moments}, for every $m\le M$, one has 
$$
\frac{1}{n} \sum_{i=1}^n  b_i^m = \frac{1}{q} \sum_{j=1}^q  c_j^m.
$$
\end{corollary} 

\proof As before, one uses the fact that the statement of Theorem \ref{moments} is invariant under parallel translations.

Let $F(x)$ be the $m$th power of the first component of vector $x$, and let $v=(t,0,\ldots,0)$. Then $F(b (x + v)) =  b^m t^m +\ldots$, where dots are terms of smaller degrees in $t$. Equating the leading terms yields  the result.
\proofend

As in the case of Steiner chains, one can say more: the value of the sum in Corollary \ref{desmom} is a symmetric polynomial of degree $m$ in two variables, the curvatures of the two parent spheres of the configuration. We do not dwell on this strengthening of Corollary \ref{desmom} here.

Theorem \ref{moments} implies Theorem \ref{Chain2}, case $n=m$ (that, in turn, implies the other cases): as a spherical design, one takes our first example with $d=1$, chooses the curvature of the circles so that they touch, and takes as the functions $F$ the real and imaginary parts of $z^m$. 

\section{In spherical geometry} \label{spgeo}

The stereographic projection from the sphere $S^2$ to the plane $\R^2$ takes circles to circles (as always, including circles of zero curvature, that is, lines) and preserves tangency. Therefore Steiner porism holds in the spherical geometry as well. In this section we describe a spherical analog of Theorem \ref{Chain1} and related topics. 

Start with the basics.
Consider the unit origin based sphere $S^d \in \R^{d+1}$. A cooriented sphere in $S^d$ is its intersection with the cooriented hyperplane  $p\cdot y = \cos\alpha$, where $p=(p_0,p_1,\ldots,p_d)$ is a unit vector, $y=(y_0,y_1,\ldots,y_d)$ are Cartesian coordinates in $\R^{d+1}$, and $\alpha$ is the radius of the sphere in the spherical metric. Denote this sphere by $C(p,\alpha)$. 
The coorientation is given by the vector $p$,  the center of $C(p,\alpha)$. The (signed) curvature of $C(p,\alpha)$ is $\cot \alpha$.

The Descartes Circle Theorem has a spherical version due to Mauldon \cite{Mau}: the curvatures of four mutually tangent cooriented circles in $S^2$ satisfy
$$
\sum_{i=1}^4 (\cot \alpha_i) = \frac{1}{2} (\sum_{i=1}^4 \cot \alpha_i)^2 -2.
$$
As in the plane, this implies that the sum of the first and  the second momenta of the curvatures remains constant in a Steiner chain of length three. 

To generalize this to Steiner chains of arbitrary length and to provide an analog of Corollary \ref{desmom}, we follow 
\cite{LMW} and use stereographic projection. 

Namely, project $S^d$ from the South Pole to the horizontal hyperplane $\R^d$ given by $y_0=0$ (when talking about stereographich projection, we always mean this one). Let $x=(x_1,\ldots,x_d)$ be Cartesian coordinates in $\R^d$. Then the image of the sphere $C(p,\alpha)$ is a sphere in $\R^d$ with  center $x$ and curvature $\bar b$ given by the formulas
\begin{equation*} 
x_i=\frac{p_i}{p_0+\cos\alpha},\ i=1,\ldots,d;\ \bar b=\frac{p_0+\cos\alpha}{\sin\alpha}=\frac{p_0}{\sin\alpha}+\cot \alpha
\end{equation*}
(p. 350 of \cite{LMW} or a direct calculation).
Inverting these formulas, we find that the stereographic projection to $S^d$ of the  sphere in $\R^d$ with center $x$ and curvature $\bar b$ has curvature given by the formula
\begin{equation} \label{curvform}
\cot \alpha = \frac{\bar b^2-1}{2\bar b} + \frac{\bar b}{2} |x|^2.
\end{equation}

Consider a sphere (with some center and some radius) in $\R^d$, and let $\{x_1,\ldots,x_n\}$ and $\{y_1,\ldots,y_q\}$ be two spherical $M$-designs. Use these points  as the centers of congruent spheres of some curvature $\bar b$.
Project these  configurations of congruent spheres stereographically from $\R^d$ to $S^d \subset \R^{d+1}$,  and let $b_i,\ i=1,\ldots,n$, and $c_j,\ j=1,\ldots,q$, be the spherical curvatures of the resulting two configurations of spheres in $S^d$. 

The next result is an analog of Corollary \ref{desmom}.

\begin{theorem} \label{spheremoments}
For every $m\le M$, one has 
$$
\frac{1}{n} \sum_{i=1}^n  b_i^m = \frac{1}{q} \sum_{j=1}^q c_j^m.
$$
\end{theorem}

\proof The argument is similar to that in the proof of Theorem \ref{moments}. 
If $x$ lies on a sphere with center $A\in \R^d$ and radius $r$, then $x=A+rv$, where $v$ is a unit vector parameterizing the sphere. It follows that $|x|^2= A^2+r^2 + 2r A\cdot v$,  a linear function on the unit sphere. Therefore the curvature (\ref{curvform}) is a linear function as well, and the proof concludes similarly to Theorem \ref{moments}. 
\proofend

Now consider a Steiner chain of length $k$ on the sphere $S^2$. We claim that there exists a rotation of the sphere such that the rotated parent circles of the Steiner chain stereographically project to concentric circles. This follows from the next lemma.

\begin{lemma} \label{concentr}
Let $A$ and $C$ be two disjoint spheres in $S^d$. There exists a rotation  of the sphere $S^d$ such that the stereographic projections of the rotated spheres $A$ and $C$ are concentric.
\end{lemma}

\proof For the proof, let us consider the sphere $S^d$ as the sphere at infinity of $d+1$-dimensional hyperbolic space in the hemisphere model. In this model, hyperbolic space is represented by the hemisphere in $\R^{d+2}$ given in Cartesian coordinates $x_1,\ldots,x_{d+2}$ by 
$$
x_1^2+\ldots +x_{d+2}^2 =1,\ x_{d+2} >0,
$$
with the metric 
$$
\frac{d x_1^2+\ldots +dx_{d+2}^2}{x_{d+2}^2}.
$$

Consider also the half-space model of hyperbolic space, given in the same coordinates by
$
x_1=1, x_{d+2} >0,
$
with the metric 
$$
\frac{d x_2^2+\ldots +dx_{d+2}^2}{x_{d+2}^2}.
$$
The stereographic projection from a point of the boundary $S^d$ takes one model to another, see \cite{Can}.

The codimension 1 spheres $A$ and $C$ in $S^d$ are boundaries of  flat (totally geodesic) hypersurfaces in the hemisphere model of hyperbolic space. Consider the geodesic perpendicular to these hypersurfaces and rotate so that its end point is the South pole of $S^d$. 

The stereographic projection takes this common perpendicular to a vertical line in the half-space model, and the two flat hypersurfaces project to two hemispheres perpendicular to this vertical line. Hence these hemispheres are concentric, and so are their boundaries, the stereographic projections of $A$ and $C$ to $\R^d$, the boundary of the half-space model.
\proofend

Lemma \ref{concentr} implies an analog of Theorem \ref{Chain1} for spherical Steiner chains.

\begin{corollary} \label{sphchain}
The first $k-1$ moments of the curvatures of the circles remain constant in the 1-parameter family of spherical Steiner chains of length $k$.
\end{corollary} 

\section{In hyperbolic geometry, briefly} \label{hypgeo}

As it happens, many a result in spherical geometry has a hyperbolic version (in formulas, one usually replaces trigonometric functions with their hyperbolic counterparts). See \cite{AP} for a general discussion of this ``analytic continuation". We briefly describe the relevant changes in the situation at hand; once again, we follow \cite{LMW}.

One uses the hyperboloid model of hyperbolic geometry: $H^d$ is represented by the upper sheet of the hyperboloid $y_0^2=y_1^2+\ldots+y_n^2+1$ in the Minkowski space $\R^{d,1}$ (again, see, e.g., \cite{Can} for various models of  hyperbolic space). The spheres in $H^d$ are the intersections of the hyperboloid with affine hyperplanes, and coorientation of a sphere is given by the coorientation of the respective hyperplane. 

The stereographic projections from ``the South Pole" $(-1,0,\ldots,0)$ to the horizontal hyperplane $y_0=0$ takes $H^d$ to 
the unit disk, and one obtains the Poincar\'e disk model of the hyperbolic space. The model is conformal, and the spheres are taken to Euclidean spheres. In particular, the Steiner porism holds in the hyperbolic plane.

The curvature of a hyperbolic sphere of radius $\alpha$ is $\coth \alpha$. For example, the hyperbolic version of the
Descartes Circle Theorem, due to Mauldon \cite{Mau}, reads:
$$
\sum_{i=1}^4 (\coth \alpha_i) = \frac{1}{2} (\sum_{i=1}^4 \coth \alpha_i)^2 +2.
$$
See \cite{LMW} for higher-dimensional generalizations. 

With this replacement of $\cot \alpha$ by $\coth \alpha$, analogs of Theorem \ref{spheremoments} and Corollary \ref{sphchain} hold in hyperbolic geometry. We do not dwell on details of proofs here.

\bigskip

{\bf Acknowledgements}. We thank A. Akopyan and J. Lagarias for their interest and encouragement.   RES and ST were  supported by NSF Research Grants DMS-1204471 and DMS-1510055, respectively. Part of this material is based upon work supported by the National Science Foundation under Grant DMS-1440140 while ST was in residence at the Mathematical Sciences Research Institute in Berkeley, California, during the Fall 2018 semester.

\end{document}